\documentclass[11pt,bezier]{article}
\usepackage{amsmath}
\usepackage{amsfonts,amsthm,amssymb}
\usepackage{amsfonts}
\usepackage{graphics}
\usepackage{amssymb}
\newtheorem{lem}{Lemma}[section]
\newtheorem{thm}[lem]{Theorem}

\theoremstyle{definition}

\begin{document}
\title{The connectivity of a bipartite graph and its bipartite complementary graph\footnote{The research is supported by NSFC (Nos. 11861066, 11531011).}}
\author{Huaping Ma, Yingzhi Tian
\footnote{Corresponding author. E-mail: mahpxj@163.com (H. Ma), tianyzhxj@163.com (Y. Tian), wly95954@163.com.}, Liyun Wu \\
{\small College of Mathematics and System Sciences, Xinjiang
University, Urumqi, Xinjiang 830046, PR China}
}

\date{}

\maketitle

\noindent{\bf Abstract } In 1956, Nordhaus and Gaddum gave lower and upper bounds on the sum and the product of the chromatic number of a graph and its complement, in terms of the order of the graph. Since then,
any bound on the sum and/or the product of an invariant in a graph $G$ and the same invariant in the complement $G^c$ of $G$ is called a Nordhaus-Gaddum type inequality or relation. The Nordhaus-Gaddum type inequalities for connectivity have been studied by several authors.
For  a  bipartite graph $G=G[X,Y]$ with bipartition ($X,Y$), its bipartite complementary graph $G^{bc}$ is a bipartite graph with $V(G^{bc})=V(G)$ and $E(G^{bc})=\{xy:\ x\in X,\ y\in Y$ and $xy \notin E(G)\}$. In this paper,  we obtain the Nordhaus-Gaddum type inequalities for connectivity of bipartite graphs and its bipartite complementary graphs. Furthermore, we prove that these inequalities are best possible.

\noindent{\bf Keywords:} Edge-connectivity; Vertex-connectivity; Bipartite graphs; Bipartite complementary graphs; Nordhaus-Gaddum type inequalities

\section{Introduction}

For graph-theoretical terminologies and notation not defined here, we follow \cite{Bondy}. In this paper, we consider finite simple graphs.  Let $G$ be a graph with vertex set $V(G)$ and edge set $E(G)$. The $order$ $n = |V(G)|$ of $G$ is the number of its vertices, while the $size$ $m=|V(G)|$ of $G$ is the number of its edges. For each vertex $v\in V(G)$, the $neighborhood$
$N(v)=N_G(v)$ of $v$ is defined as the set of all vertices adjacent to $v$, and $d(v)=|N(v)|$ is the $degree$ of $v$. We denote by $\delta(G)$
the $minimum$ $degree$ and by $\Delta(G)$ the $maximum$ $degree$ of $G$.
The $vertex$-$connectivity$ $\kappa(G)$ of the graph $G$ is the minimum size of a  vertex set $S$  such that $G-S$ is disconnected or has only one vertex; the $edge$-$connectivity$ $\kappa'(G)$ of the graph $G$ is  the minimum size of an edge set $F$  such that $G-F$ is disconnected. The graph $G$ is said to be $k$-$vertex$-$connected$ ($k$-$edge$-$connected$) if $\kappa(G)\geq k$ ($\kappa'(G)\geq k$).
For $\kappa(G)\leq\kappa'(G)\leq\delta(G)$, a graph $G$ with $\kappa(G)=\delta(G)$
is  called $maximally\ vertex$-$connected$;  a graph $G$ with $\kappa'(G)=\delta(G)$ is  called $maximally\ edge$-$connected$.
The complement $G^c$ of $G$ is the graph defined
on the vertex set $V(G)$ of $G$, where an edge $uv$ belongs to $G^c$ if and only if it does not belong to $G$. The $floor$ of a real number $x$, denoted by $\lfloor x\rfloor$, is the greatest integer not larger than $x$; the $ceil$ of a real number $x$, denoted by $\lceil x\rceil$, is the least integer greater than or equal to $x$.

A $coloring$ of $G$ is an assignment of colors to the vertices of $G$ such that two adjacent vertices have different colors. The $chromatic$ $number$ of $G$, denoted by $\chi(G)$, is the minimum number of colors in a coloring of $G$.
In 1956, Nordhaus and Gaddum \cite{Nordhaus} gave lower and upper bounds on the sum and on the product of $\chi(G)$ and $\chi(G^c)$ in terms of the order $n$ of $G$. Since then, any bound on the sum and/or the product
of an invariant in a graph $G$ and the same invariant in the complement $G^c$ of $G$ is called a Nordhaus-Gaddum type inequality
or relation. For many  graph invariants, relations of a similar type have been proposed, see the survey \cite{Aouchiche}.

The original relations presented by Nordhaus and Gaddum  are as follows.

\begin{thm}
(\cite{Nordhaus}) If $G$ is a graph of order $n$, then
$$2\sqrt n\leq \chi(G)+\chi(G^c)\leq n+1\
and
\ n\leq \chi(G)\chi(G^c)\leq \frac{(n+1)^2}{4}.$$
Furthermore, these bounds are best possible for infinitely many values of n.
\end{thm}

In 1970, Alavi and Mitchem \cite{Alavi} proved all the Nordhaus-Gaddum type
inequalities for both vertex-connectivity and edge-connectivity.

\begin{thm}
(\cite{Alavi}) For any graph $G$ on $n\ (\geq2)$ vertices,
$$1\leq \kappa'(G)+\kappa'(G^c)\leq n-1,$$
and
\[
0\leq \kappa'(G)\kappa'(G^c)\leq
\left\{
\begin{array}{ll}
\lfloor\frac{n-1}{2}\rfloor\lceil\frac{n-1}{2}\rceil, &n\equiv0,1,2\ (\ mod\ 4\ ),\\
\frac{n-3}{2}\frac{n+1}{2},& n\equiv3\ (\ mod\ 4\ ).
\end{array}
\right.
\]
Furthermore, these bounds are sharp.
\end{thm}

Due to the reason that many of the sharp bounds in the classical Nordhaus-Gaddum type inequalities are attained by one of the graphs $G$ and $G^c$ being very dense. Achuthan et al. \cite{Achuthan} in 1990 considered the Nordhaus-Gaddum problem by restricting $G$ to the graphs with given size $m$.

\begin{thm}(\cite{Achuthan})
If $G$ is a graph on $n$ vertices and $m$ edges, where $m\leq\frac{n(n-1)}{4}$, then
\[
max\{1,n-1-m\}\leq\kappa'(G)+\kappa'(G^c)\leq
\left\{
\begin{array}{ll}
n-3, &n+1\leq 2m\leq 2n-4,\\
n-2,& 2\leq 2m \leq n,\ or\ m=n-1,\\
    &or\ 2m\not\equiv 0\ (\ mod\ n\ )\ and\ m\geq n,\\
n-1,&otherwise,\\
\end{array}
\right.
\]
and
\[
0\leq\kappa'(G)\kappa'(G^c)\leq
\left\{
\begin{array}{ll}
0, &m\leq n-2,\\
\frac{2m}{n}(n-1-\frac{2m}{n}),& 2m \equiv 0(mod\ n\ )\ and\ m\geq n,\\
\lfloor\frac{2m}{n}\rfloor(n-2-\lfloor\frac{2m}{n}\rfloor),&otherwise.\\
\end{array}
\right.
\]
Furthermore, these bounds are sharp for all $n$ and $m$.
\end{thm}

In 2008, Hellwig and Volkmann \cite{Hellwig} proved that if $G$ and $G^c$ are connected graphs, then $\kappa'(G)+\kappa'(G^c)\geq min\{\delta(G),\delta(G^c)\}+1$ and $\kappa(G)+\kappa(G^c)\geq min\{\delta(G),\delta(G^c)\}+1$. Moreover, these inequalities are best possible.

For a bipartite graph $G=G[X,Y]$  with  bipartition ($X$,$Y$),  the $bipartite$ $complementary$ $graph$ $G^{bc}$ of  $G$ is defined  to be the bipartite graph  with $V(G^{bc})=V(G)$ and
$E(G^{bc})=\{xy:\ x\in X,\ y\in Y$ and $xy \notin E(G)\}$.
A $complete$ $bipartite$ $graph$  is a special kind of bipartite graph  where every vertex of $X$ is connected to every vertex of $Y$. A complete bipartite graph with partitions of size $|X|=r$  and $|Y|=s$, is denoted by $K_{r,s}$.

Motivated by the Nordhaus-Gaddum type inequalities for the connectivity of a graph and its complement in Theorems 1.2 and 1.3, we consider similar results for the bipartite graphs and its bipartite complementary graphs in this paper. The graphs in Section 2 will be used to prove the sharpness of the inequalities. In section 3, similar results as Theorem 1.2 are obtained for the bipartite graphs and its bipartite complementary graphs. We give similar results as Theorem 1.3 for the bipartite graphs and its bipartite complementary graphs in the last section.

\section{Preliminaries}

A graph $G$ is said to be $vertex$-$transitive$
($edge$-$transitive$) if $Aut(G)$ acts transitively on $V(G)$
($E(G)$), that is, for any two vertices $u$ and $v$ (two edges $e_1$
and $e_2$) in $G$, there is an automorphism $\alpha$ of $G$ such
that $v=\alpha(u)$ ($e_2=\alpha(e_1)$). A bipartite graph $G=G[X,Y]$ with bipartition $X\cup Y$ is called $half\ vertex$-$transitive$  if $Aut(G)$ acts
transitively both on $X$ and $Y$. Let $Z_n$ be the cyclic group of integers modulo $n$.

{\flushleft\textbf{Definition 2.1.}} For a group $A$, let $S$ be a
subset of $A$ such that $1_A\notin S\ and\  S^{-1}=S $, the Cayley
graph $Cay(A,S)$ is a graph with vertex set $A$ and edge set
$\{\{g,sg\}:g\in A, s\in S\}$.

It is well known that Cayley graphs are vertex-transitive. For
studying semi-symmetric graphs, which are regular edge-transitive
but not vertex-transitive, Xu defined the Bi-Cayley graphs.

{\flushleft\textbf{Definition 2.2.}} (\cite{Xu}) For a group $A$, let
$S$ be a subset of $A$, the $Bi$-$Cayley$ $graph$ $BC(A,S)$ is a bipartite graph with vertex set $A \times \{0,1\}$  and edge set
$\{\{(g,0),(sg,1)\}:g\in A, s\in S\}$.

The following lemma is easily obtained by Definition 2.2.

\begin{lem}
Let $A$ be a group and let $S$ be a subset of $A$. Then the bipartite complement $G^{bc}$ of the Bi-Cayley graph $G=BC(A,S)$ is also a Bi-Cayley graph. Moreover, $G^{bc}\cong BC(A,A\setminus S)$.
\end{lem}

For $a\in A$, the translation $MR(a)$ defined by $(x,i)\rightarrow (xa,i) \ (i=0,1)$ is an automorphism of the Bi-Cayley graph $G=BC(A,S)$. Since all of these automorphisms form a subgroup $MR(A)$ of $Aut(G)$, which acts transitively both on $A\times\{0\}$ and $A\times\{1\}$, thus $G=BC(A,S)$ is half vertex-transitive.

For vertex-transitive graphs, the following result are well-known.

\begin{thm}(\cite{Mader})
All connected vertex-transitive graphs are maximally edge-connected.
\end{thm}

For half vertex-transitive graphs, Liang and Meng \cite{Liang}
proved the following result.

\begin{thm} (\cite{Liang})
Every connected half vertex-transitive graph $G$ is maximally vertex-connected, and thus is maximally edge-connected.
\end{thm}

Combining Lemma 2.1 with Theorem 2.3, we have the following lemma.

\begin{lem}
Let $A$ be a group and let $S$ be a subset of $A$. If both the Bi-Cayley graph $G=BC(A,S)$ and its bipartite complement $G^{bc}$ are connected, then both $G$ and $G^{bc}$ are maximally vertex-connected, and thus are maximally edge-connected.
\end{lem}

\begin{lem}
Let $G$ be a $k$-edge-connected graph.  If  $G'$ is a graph obtained from $G$ by adding a new vertex $v$ and at least $k$ edges between $v$ and $V(G)$, then $G'$ is also $k$-edge-connected.
\end{lem}

\noindent{\bf Proof.} For any minimum edge-cut $F$ of $G'$, either $F\cap E(G)$ is an edge-cut of $G$ or $F$ is the set of edges incident with $v$,
thus $G'$ is $k$-edge-connected. $\Box$

\section{The connectivity of a bipartite graph and its bipartite complement}

We first give bounds on the sum and the product for the minimum degrees of a bipartite graph and its bipartite complement.

\begin{lem}
Let $G=G[X,Y]$ be a bipartite graph on $n$ vertices. Assume $|X|=r$, $|Y|=s$ and $r\leq s$. Then

$$0\leq\delta(G)+\delta(G^{bc})\leq r,  \eqno(1)$$
and
$$0\leq\delta(G)\delta(G^{bc})\leq\lceil\frac{r}{2}\rceil
\lfloor\frac{r}{2}\rfloor. \eqno(2) $$
\end{lem}

\noindent{\bf Proof.} It is clear that $0\leq\delta(G)+\delta(G^{bc})\leq r$ and $0\leq\delta(G)\delta(G^{bc})$. If the sum of two numbers is $r$, then $\delta(G)\delta(G^{bc})$ is maximum when each of them is $\frac{r}{2}$.

If $r$ is even, then $\frac{r}{2}$ is an integer and
$$\delta(G)\delta(G^{bc})\leq(\frac{\delta(G)+\delta(G^{bc})}{2})^2\leq
(\frac{r}{2})^2,$$
so  $\delta(G)\delta(G^{bc})\leq\lceil\frac{r}{2}\rceil
\lfloor\frac{r}{2}\rfloor$ is obtained. On the other hand, if $r$ is odd, then $\frac{r}{2}$ is not an integer. However,
$$(\frac{r}{2})^2=(\frac{r+1}{2})(\frac{r-1}{2})+\frac{1}{4},$$
where $\frac{r+1}{2}$ and $\frac{r-1}{2}$ are integers. So
$\delta(G)\delta(G^{bc})\leq\lceil\frac{r}{2}\rceil
\lfloor\frac{r}{2}\rfloor$ holds.  $\Box$

Now we are ready to give the Nordhaus-Gaddum type inequalities for the edge-connectivity of a bipartite graph and it bipartite complement.

\begin{thm}
Let $G=G[X,Y]$ be a bipartite graph on $n$ vertices. Assume $|X|=r$, $|Y|=s$ and $r\leq s$. Then

$$0\leq \kappa'(G)+\kappa'(G^{bc})\leq r, \eqno(3)$$
and
$$0\leq \kappa'(G)\kappa'(G^{bc})\leq\lceil\frac{r}{2}\rceil\lfloor\frac{r}{2}\rfloor.
\eqno(4)$$
Moreover, the lower bound of (3) is sharp for all $n\geq3$, the upper bound of (3) and the lower bound of (4) are sharp for all $n\geq2$, and the upper bound of (4) is sharp for all $r\geq4$.
\end{thm}

\noindent{\bf Proof.} The lower bounds of (3) and (4) are immediate.  Since $\kappa'(G)\leq\delta(G)$, we have $\kappa'(G)+\kappa'(G^{bc})\leq\delta(G)+\delta(G^{bc})\leq r$ and $\kappa'(G)\kappa'(G^{bc})\leq\delta(G)\delta(G^{bc})
\leq\lceil\frac{r}{2}\rceil\lfloor\frac{r}{2}\rfloor$ by Lemma 3.1. Thus the upper bounds of (3) and (4) hold.

Let $G_1=G_1[X,Y]$ be a bipartite graph on $n$ ($n\geq3$) vertices such that ($i$) there is a vertex $y_1$ in $Y$ adjacent to all vertices in $X$ and ($ii$) there is another vertex $y_2$ in $Y$ not adjacent to any vertices in $X$. Then both $G_1$ and $G_1^{bc}$ are disconnected. Therefore, the lower bound of (3) is sharp for all $n\geq3$.

The complete bipartite graph on $n$ vertices shows that the upper bound of (3) and the lower bound of (4) are best possible for all $n\geq2$.

To prove the sharpness of the upper bound of (4), we construct a bipartite graph $G_2$ from $BC(Z_r,\{0,1,\cdots,\lfloor\frac{r}{2}\rfloor-1\})$ by adding $s-r$ vertices $y_{r+1},\cdots,y_{s}$, where $y_i$ is adjacent to exact $\lfloor\frac{r}{2}\rfloor$ vertices in  $Z_r\times\{0\}$ for $i=r+1,\cdots,s$. Then $G_2^{bc}$ is obtained from $BC(Z_r,\{\lfloor\frac{r}{2}\rfloor,\lfloor\frac{r}{2}\rfloor+1,\cdots,r-1\})$ by adding $s-r$ vertices $y_{r+1},\cdots,y_{s}$, where $y_i$ is adjacent to each vertex in $Z_r\times\{0\}\setminus N_{G_2}(y_i)$ for $i=r+1,\cdots,s$. Since $r\geq4$, both $G_2$ and $G_2^{bc}$ are connected. By Lemmas 2.4 and 2.5, $\kappa'(G_2)=\lfloor\frac{r}{2}\rfloor$ and $\kappa'(G_2^{bc})=\lceil\frac{r}{2}\rceil$. Thus the upper bound of (4) is best possible for all $r\geq4$.
$\Box$

Since Lemma 2.5 is also true for vertex-connectivity, by a similar argument as Theorem 3.2, we present the following theorem without proof.

\begin{thm}
Let $G=G[X,Y]$ be a bipartite graph on $n$ vertices. Assume $|X|=r$, $|Y|=s$ and $r\leq s$. Then
$$0\leq \kappa(G)+\kappa(G^{bc})\leq r,$$
and
$$0\leq \kappa(G)\kappa(G^{bc})\leq\lceil\frac{r}{2}\rceil\lfloor\frac{r}{2}\rfloor.$$
Furthermore, these bounds are best possible.
\end{thm}

\section{The connectivity of a bipartite graph with given size $m$ and its bipartite complement}

In this section, let $G=G[X,Y]$ be a bipartite graph with $n$ vertices and $m$ edges, where $X=\{x_1,\cdots,x_r\}$ and $Y=\{y_1,\cdots,y_s\}$.  Without loss of generality, we assume that $r\leq s$ and $m\leq\frac{1}{2}rs$.

\begin{thm}
Let $G=G[X,Y]$ be a bipartite graph with $n$ vertices and $m$ edges. Assume $|X|=r$, $|Y|=s$, $r\leq s$ and $m\leq\frac{1}{2}rs$. Then
$$max \{0, r-m\}\leq\kappa'(G)+\kappa'(G^{bc})\leq N(n,m), \eqno(5)$$
where
\[
N(n,m)=
\left\{
\begin{array}{ll}
r-2, &s+1\leq m\leq n-2,\\
r-1, &1\leq m\leq s,\ or\ m=n-1\ and\ r\geq2,\ or\ m\not\equiv  0(mod\ s)\ and\ m\geq n, \\
r,& otherwise.
\end{array}
\right.
\]
Furthermore, these bounds are best possible for all $n$ and $m$.
\end{thm}

\noindent{\bf Proof.} Clearly,  $\kappa'(G)+\kappa'(G^{bc})\geq 0$. The lower bound of (5) holds when $m \geq r$. So let $m< r$, we shall show that $\kappa'(G)+\kappa'(G^{bc})\geq r-m$.
Since $m< r$, we know $G$ is disconnected and $\kappa'(G)=0$. Thus we need to show that $\kappa'(G^{bc})\geq r-m$. We note that $G^{bc}$ is isomorphic to a bipartite graph obtained from the complete bipartite graph $K_{r,s}$ by deleting $m$ edges. By deleting one edge from a graph, its edge-connectivity decreases at most one. Therefore $\kappa'(G^{bc})\geq \kappa'(K_{r,s})-m=r-m$, and the lower bound of (5) holds.

By Theorem 3.2, $\kappa'(G)+\kappa'(G^{bc})\leq r$. If $1\leq m\leq n-2$, then $\kappa'(G)=0$ and $\kappa'(G^{bc})\leq\delta(G^{bc})\leq r-1$. Thus $\kappa'(G)+\kappa'(G^{bc})\leq r-1$. If $m=n-1$ and $r\geq2$, then there is vertex $y\in Y$ such that $d_G(y)\geq2$, which implies $d_{G^{bc}}(y)\leq r-2$. Thus $\kappa'(G)+\kappa'(G^{bc})\leq 1+ r-2=r-1$. If $m\neq 0(mod\ s)$ and $m\geq n$, let $m=ds+l$, where $1\leq l\leq s-1$. Then there is a vertex $y_1\in Y$ such that $d_G(y_1)\leq d$ and there is a vertex $y_2\in Y$ such that $d_G(y_2)\geq d+1$. Thus $\kappa'(G)+\kappa'(G^{bc})\leq d_G(y_1)+d_{G^{bc}}(y_2)\leq d+(r-d-1)=r-1$. If $s+1\leq m\leq n-2$, then there is a vertex $y_3\in Y$ such that $d_G(y_3)\geq 2$. Thus $\kappa'(G)+\kappa'(G^{bc})\leq 0+d_{G^{bc}}(y_3)\leq 0+(r-2)=r-2$. This implies that the upper bounds of (5) are true.

In the following, we first prove the sharpness of the lower bound. If $m< r$, define $G_1$ to be the bipartite graph with bipartition $X\cup Y$ and edge set $\{x_iy_1:1\leq i\leq m\}$. Then $\kappa'(G_1)=0$ and $\kappa'(G_1^{bc})=r-m$. If $m\geq r$, define $G_2$ to be a bipartite graph with bipartition $X\cup Y$ and define the edge set as follows: ($i$) $y_1$ is adjacent to each vertex in $X$; ($ii$) the remaining $m-r$ edges are connected between $X$ and $\{y_3,\cdots,y_s\}$ (this can be done by the assumption that $m\leq\frac{1}{2}rs$). Since $y_1$ is adjacent to each vertex in $X$ and $y_2$ is an isolated vertex in $G_2$, we have $\kappa'(G_2)=0$ and $\kappa'(G_2^{bc})=0$. Thus the lower bound of (5) are best possible for all $n$ and $m$.

Now we are ready to prove the sharpness of the upper bounds. If $m=0$, then $G$ is an empty graph and $G^{bc}$ is a complete bipartite graph, which implies that $\kappa'(G)=0$, $\kappa'(G^{bc})=r$ and $\kappa'(G)+\kappa'(G^{bc})=r$. If $1\leq m\leq s$, let $G_3$ be a bipartite graph with bipartition $X\cup Y$ and edge set $\{x_iy_i:1\leq i\leq r\}\cup\{x_1y_j:r+1\leq j\leq m\}$. Then $\kappa'(G_3)=0$, $\kappa'(G_3^{bc})=r-1$ and $\kappa'(G_3)+\kappa'(G_3^{bc})=r-1$.
If $s+1\leq m\leq n-2$, let $G_4$ be a bipartite graph with bipartition $X\cup Y$ and edge set $\{x_iy_i:1\leq i\leq r\}\cup\{x_iy_{i+1}:1\leq i\leq m-s\}\cup\{x_1y_i:r+1\leq i\leq s\}$, where $m-s\leq n-2-s=r-2$.
Then $\kappa'(G_4)=0$, $\kappa'(G_4^{bc})=r-2$ and $\kappa'(G_4)+\kappa'(G_4^{bc})=r-2$.
If $m=n-1$ and $r=1$, then $G$ is isomorphic to $K_{1,n-1}$ and
$G^{bc}$ is an empty graph, which implies that $\kappa'(G)=r=1$, $\kappa'(G^{bc})=0$ and $\kappa'(G)+\kappa'(G^{bc})=r$.
If $m=n-1$ and $r\geq2$, let $G_5$ be a bipartite graph with bipartition $X\cup Y$ and edge set $\{x_iy_i:1\leq i\leq r\}\cup\{x_iy_{i+1}:1\leq i\leq r-1\}\cup\{x_ry_j:r+1\leq j\leq s\}$. Then $\kappa'(G_5)=1$, $\kappa'(G_5^{bc})=r-2$ and $\kappa'(G_5)+\kappa'(G_5^{bc})=r-1$.

If $m\equiv 0(mod\ s)$ and $m\geq n$, let $m=ds$. We construct a bipartite graph $G_6$ from $BC(Z_r,\{0,1,\cdots,d-1\})$ by adding $s-r$ vertices $y_{r+1},\cdots,y_{s}$, where $y_i$ is adjacent to exact $d$ vertices of  $Z_r\times\{0\}$ for $i=r+1,\cdots,s$. Then $G_6^{bc}$ is obtained from $BC(Z_r,\{d, d+1,\cdots,r-1\})$ by adding $s-r$ vertices $y_{r+1},\cdots,y_{s}$, where $y_i$ is adjacent to each vertex in $Z_r\times\{0\}\setminus N_{G_6}(y_i)$ for $i=r+1,\cdots,s$. Since $d\geq2$ and $m\leq\frac{1}{2}rs$, both $G_6$ and $G_6^{bc}$ are connected. By Lemmas 2.4 and 2.5, $\kappa'(G_6)=d$ and $\kappa'(G_6^{bc})=r-d$. Thus $\kappa'(G_6)+\kappa'(G_6^{bc})=r$.

If $m\not\equiv 0(mod\ s)$ and $m\geq n$, let $m=ds+l$, where $1\leq l\leq s-1$. We construct a bipartite graph $G_7$ from $G_6$ by adding edges $x_1y_i$, $d\leq i\leq d+l-1$. Then $\kappa'(G_7)=d$ and $\kappa'(G_7^{bc})=r-d-1$. Thus $\kappa'(G_7)+\kappa'(G_7^{bc})=r-1$.

The theorem is thus established.
$\Box$

\begin{thm}
Let $G=G[X,Y]$ be a bipartite graph with $n$ vertices and $m$ edges. Assume $|X|=r$, $|Y|=s$, $r\leq s$ and $m\leq\frac{1}{2}rs$. Then
$$0\leq \kappa'(G)\kappa'(G^{bc})\leq M(n,m), \eqno(6)$$
where
\[
M(n,m)=
\left\{
\begin{array}{ll}
0, &m\leq n-2,\ or\ m=n-1\ and\ r=1,\\
\frac{m}{s}(r-\frac{m}{s}), & m\ \equiv 0\  (mod \ s) \ and\ m\ \geq\ n, \\
\lfloor\frac{m}{s}\rfloor(r-1-\lfloor\frac{m}{s}\rfloor),& otherwise.
\end{array}
\right.
\]
Furthermore, these bounds are best possible for all $n$ and $m$.
\end{thm}

\noindent{\bf Proof.} Clearly, $\kappa'(G)\kappa'(G^{bc})\geq0$.

If $m\leq n-2$, then $G$ is disconnected and $\kappa'(G)=0$. Thus $\kappa'(G)\kappa'(G^{bc})=0$. If $m=n-1$ and $r=1$, then $G$ is isomorphic to $K_{1,n-1}$ and  $G^{bc}$ is an empty graph. Therefore, $\kappa'(G)\kappa'(G^{bc})=0$. If $m=n-1$ and $r\geq2$, then there is vertex $y\in Y$ such that $d_G(y)\geq2$, which implies $d_{G^{bc}}(y)\leq r-2$. By $\kappa'(G)\leq1$, we obtain that $\kappa'(G)\kappa'(G^{bc})\leq 1\times(r-2)=\lfloor\frac{m}{s}\rfloor(r-1-\lfloor\frac{m}{s}\rfloor)$.

If $m\equiv 0\  (mod \ s)$ and $m\geq n$, say $m=sd$, then there is a vertex $y_1\in Y$ with $d_G(y_1)\leq d$ and  there is  a vertex $y_2\in Y$ with $d_G(y_2)\geq d$. Thus
$\delta(G)\leq d$, $\delta(G^{bc})\leq r-d$ and
$\kappa'(G)\kappa'(G^{bc})\leq \frac{m}{s}(r-\frac{m}{s})$.

If $m\not\equiv 0\  (mod \ s)$ and $m\geq n$, say $m=sd+l$, where $1\leq l \leq s-1$, then there is a vertex $y_3\in Y$ with $d_G(y_3)\leq d$ and  there is  a vertex $y_4\in Y$ with $d_G(y_4)\geq d+1$. Thus
$\delta(G)\leq d$, $\delta(G^{bc})\leq r-d-1$ and $\kappa'(G)\kappa'(G^{bc})\leq \lfloor\frac{m}{s}\rfloor(r-1-\lfloor\frac{m}{s}\rfloor)$.

The bipartite graphs $G_1$ and $G_2$ constructed in the proof of Theorem 4.1 show that the lower bound of (6) are best possible for all $n$ and $m$.

In the following, we will prove the sharpness of the upper bounds.  If $m\leq n-2$, then $G$ is disconnected. Thus $\kappa'(G)\kappa'(G^{bc})=0$.
If $m=n-1$ and $r=1$, then $G$ is isomorphic to $K_{1,n-1}$ and
$G^{bc}$ is an empty graph, which implies that $\kappa'(G)=r=1$, $\kappa'(G^{bc})=0$ and $\kappa'(G)\kappa'(G^{bc})=0$.
If $m=n-1$ and $r\geq2$, the bipartite graph $G_5$ constructed in the proof of Theorem 4.1 shows that $\kappa'(G)\kappa'(G^{bc})= 1\times(r-2)=\lfloor\frac{m}{s}\rfloor(r-1-\lfloor\frac{m}{s}\rfloor)$.

If $m\equiv 0(mod\ s)$ and $m\geq n$, then the bipartite graph $G_6$ constructed in the proof of Theorem 4.1 shows that $\kappa'(G)\kappa'(G^{bc})=
\frac{m}{s}(r-\frac{m}{s})$.

If $m\not\equiv 0(mod\ s)$ and $m\geq n$, then the bipartite graph $G_7$ constructed in the proof of Theorem 4.1 shows that $\kappa'(G)\kappa'(G^{bc})=
\lfloor\frac{m}{s}\rfloor(r-1-\lfloor\frac{m}{s}\rfloor)$.

The theorem is thus established.
$\Box$

Since Lemma 2.5 is also true for vertex-connectivity, by similar arguments as Theorems 4.1 and 4.2, we have the following Nordhaus-Gaddum type inequalities for vertex-connectivity of the bipartite graphs with given size $m$ and its bipartite complementary graphs.

\begin{thm}
Let $G=G[X,Y]$ be a bipartite graph with $n$ vertices and $m$ edges. Assume $|X|=r$, $|Y|=s$, $r\leq s$ and $m\leq\frac{1}{2}rs$. Then
$$max \{0, r-m\}\leq\kappa(G)+\kappa(G^{bc})\leq N(n,m),$$
and
$$0\leq \kappa(G)\kappa(G^{bc})\leq M(n,m),$$
where
\[
N(n,m)=
\left\{
\begin{array}{ll}
r-2, &s+1\leq m\leq n-2,\\
r-1, &1\leq m\leq s,\ or\ m=n-1\ and\ r\geq2,\ or\ m\not\equiv  0(mod\ s)\ and\ m\geq n, \\
r,& otherwise,
\end{array}
\right.
\]
and
\[
M(n,m)=
\left\{
\begin{array}{ll}
0, &m\leq n-2,\ or\ m=n-1\ and\ r=1,\\
\frac{m}{s}(r-\frac{m}{s}), & m\ \equiv 0\  (mod \ s) \ and\ m\ \geq\ n, \\
\lfloor\frac{m}{s}\rfloor(r-1-\lfloor\frac{m}{s}\rfloor),& otherwise.
\end{array}
\right.
\]
Furthermore, these bounds are best possible.
\end{thm}

\end{document}